\documentclass[final]{siamltex}

\usepackage[vlined,algosection,ruled]{algorithm2e}
\usepackage{tikz}
\usepackage{tikz-qtree}
\usepackage{amsmath,amssymb}
\usepackage[normalem]{ulem} 
\usepackage{verbatim} 
\usepackage{url}
\usepackage{pgfplots}

\makeatletter
\renewcommand*\env@matrix[1][\arraystretch]{%
  \edef\arraystretch{#1}%
  \hskip -\arraycolsep
  \let\@ifnextchar\new@ifnextchar
  \array{*\c@MaxMatrixCols c}}
\makeatother

\pgfplotsset{compat=newest}

\pgfplotscreateplotcyclelist{journal}{%
blue,every mark/.append style={fill=blue!80!black},mark=square*\\%
black,dashed\\%
red,densely dashed,every mark/.append style={solid,fill=red!80!black},mark=o\\%
black,dashed\\%
brown!80!black,every mark/.append style={fill=brown!80!black},mark=diamond\\%
black,dashed\\%
green,every mark/.append style={fill=green!80!black},mark=diamond*\\%
black,dashed\\%
cyan,densely dashed,every mark/.append style={solid,fill=cyan!80!black},mark=*\\%
black,dashed\\%
magenta,densely dashed,every mark/.append style={solid,fill=magenta!80!black},mark=square*\\%
black,dashed\\%
}

\title{A Parallel Butterfly Algorithm}

\author{Jack Poulson\thanks{Department of Mathematics, Stanford University, 450 Serra Mall, Stanford, CA 94305 (poulson@stanford.edu).} 
\and Laurent Demanet\thanks{Department of Mathematics, Massachusetts Institute of Technology, 77 Massachusetts Avenue, Cambridge, MA 02139 (laurent@math.mit.edu)}
\and Nicholas Maxwell\thanks{Department of Mathematics, University of Houston, 651 PGH, Houston, TX 77004 (nmaxwell@math.uh.edu)} 
\and Lexing Ying\thanks{Department of Mathematics and ICME, Stanford University, 450 Serra Mall, Stanford, CA 94305 (lexing@math.stanford.edu).} 
}

\date{\today}

\begin{document}

\maketitle

\begin{abstract}
The {\em butterfly algorithm} is a fast algorithm which approximately 
evaluates a discrete analogue of the integral transform 
$\int_{\mathbb{R}^d} K(x,y) g(y) dy$ at large numbers of target points
when the kernel, $K(x,y)$, is approximately low-rank when restricted to 
subdomains satisfying a certain simple geometric condition.
In $d$ dimensions with $O(N^d)$ quasi-uniformly distributed source and target points, 
when each appropriate submatrix of $K$ is approximately rank-$r$, the running time of 
the algorithm is at most $O(r^2 N^d \log N)$.
A parallelization of the butterfly algorithm is introduced which, assuming a 
message latency of $\alpha$ and per-process inverse bandwidth of $\beta$, 
executes in at most
$O\left(r^2 \frac{N^d}{p} \log N + 
  \left(\beta r\frac{N^d}{p}+\alpha\right)\log p\right)$ 
time using $p$ processes.
This parallel algorithm was then instantiated in the form of the open-source 
\texttt{DistButterfly} library for the special case where 
$K(x,y)=\exp(i \Phi(x,y))$, where $\Phi(x,y)$ is a black-box, sufficiently 
smooth, real-valued {\em phase function}.
Experiments on Blue Gene/Q demonstrate impressive strong-scaling results for 
important classes of phase functions. Using quasi-uniform sources, 
hyperbolic Radon transforms and an analogue of a 3D generalized Radon transform 
were respectively observed to strong-scale from 1-node/16-cores up to 
1024-nodes/16,384-cores with greater than 90\% and
82\% efficiency, respectively. 
\end{abstract}



\pagestyle{myheadings}
\thispagestyle{plain}
\markboth{J. POULSON ET AL.}{PARALLEL BUTTERFLY ALGORITHM}

\section{Introduction}
The {\em butterfly algorithm}~\cite{MichielssenBoag-multilevel,ONeil-special,Ying-SFT,CDY-butterfly} 
provides an efficient means of (approximately) applying any integral operator
\[
  (\mathcal K g)(x) = \int_Y K(x,y) g(y) dy
\]
whose {\em kernel}, $K : X \times Y \rightarrow \mathbb{C}$, satisfies the 
condition that, given any {\em source box}, $B \subset Y$, and 
{\em target box}, $A \subset X$, 
such that the product of their diameters is less than some fixed constant, 
say $D$, the restriction of $K$ to $A \times B$, henceforth denoted 
$K|_{A \times B}$, is approximately low-rank in a point-wise sense.
More precisely, there exists a numerical rank, $r(\epsilon)$, which depends at most 
polylogarithmically on $1/\epsilon$, such that, for any subdomain $A \times B$ satisfying 
$\text{diam}(A)\text{diam}(B) \le D$, there exists a rank-$r$ separable approximation 
\[
  \left| K(x,y) - \sum_{t=0}^{r-1} u_t^{AB}(x) v_t^{AB}(y) \right| 
  \le \epsilon,\;\;\;\forall\, x \in A,\, y \in B.
\]

If the source function, $g : Y \rightarrow \mathbb{C}$, is lumped into some 
finite set of points, $I_Y \subset Y$, the low-rank separable representation 
implies an approximation
\[
  \left| f^{AB}(x) - 
         \sum_{t=0}^{r-1} u_t^{AB}(x) \delta_t^{AB} \right| \le 
  \epsilon\, \| g_B \|_1,\;\;\;\forall\, x \in A,
\]
where 
\begin{equation}
  f^{AB}(x) \equiv \sum_{y \in I_Y \cap B} K(x,y) g(y),\;\;\;x \in A,
\end{equation}
represents the potential generated in box $A$ due to the sources in box $B$, 
$\| g_B \|_1$ is the $\mathcal{\ell}_1$ norm of $g$ over its support in 
$B$, and each {\em expansion weight}, $\delta_t^{AB}$, could simply be chosen
as
\[
  \delta_t^{AB} = \sum_{y \in I_Y \cap B} v_t^{AB}(y) g(y).
\]

The first step of the butterfly algorithm is to partition the source
domain, $Y$, into a collection of boxes which are sufficiently small such that
the product of each of their diameters with that of the entire target domain, 
$X$, is less than or equal to $D$. Then, for each box $B$ in this 
initial partitioning of the source domain, the expansion weights, 
$\{\delta_t^{XB}\}_{t=0}^{r-1}$, for approximating the potential over $X$
 generated by the sources in $B$, $f^{XB}$, can be cheaply computed as
$\delta_t^{XB} := \sum_{y \in I_Y \cap B} v_t^{XB}(y) g(y)$.
In $d$ dimensions, if there are $N^d$ such source boxes, each containing $O(1)$ 
sources, then this initialization step only requires $O(r N^d)$ work.

Upon completion of the butterfly algorithm, we will have access to a much more
useful set of expansion weights, those of $\{f^{AY}\}_A$, where each box $A$ 
is a member of a sufficiently fine partitioning of the target domain such that 
the product of its diameter with that of the entire source domain is bounded
by $D$.
Then, given any $x \in X$, there exists a box $A \ni x$ within which we may 
cheaply evaluate the approximate solution
\[
  f(x) = f^{AY}(x) \approx \sum_{t=0}^{r-1} u_t^{AY}(x) \delta_t^{AY}.
\]
If we are interested in evaluating the solution at $N^d$ target
points, then the final evaluation phase clearly requires $O(r N^d)$ work.

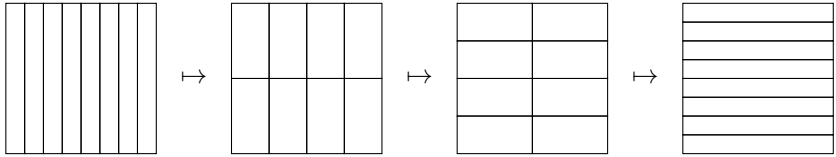
\begin{figure}
\begin{center}
\begin{tikzpicture}[scale=2]

\draw[white] (3.25,0.75) circle (0.3);

\foreach \j in {0,...,7}{
  \draw[black] (0.125*\j,0) rectangle (0.125*\j+0.125,1);
}

\draw (1.25,0.5) node {$\mapsto$};

\foreach \i in {0,1}{
  \foreach \j in {0,...,3}{
    \draw[black] (1.5+0.25*\j,0.5*\i) rectangle (1.5+0.25*\j+0.25,0.5*\i+0.5);
  }
}

\draw (2.75,0.5) node {$\mapsto$};

\foreach \i in {0,...,3}{
  \foreach \j in {0,1}{
    \draw[black] (3+0.5*\j,0.25*\i) rectangle (3+0.5*\j+0.5,0.25*\i+0.25);
  }
}

\draw (4.25,0.5) node {$\mapsto$};

\foreach \i in {0,...,7}{
  \draw[black] (4.5+0,0.125*\i) rectangle (4.5+1,0.125*\i+0.125);
}

\end{tikzpicture}
\end{center}
\caption{The successive partitions of the product space $X \times Y$ during a 
1D butterfly algorithm with $N=8$. This figure has an additional matrix-centric 
interpretation: each rectangle is an approximately low-rank submatrix of the 
discrete $1D$ integral operator.}
\label{fig:1D}
\end{figure}

The vast majority of the work of the butterfly algorithm lies in 
the translation of the expansion weights used to approximate the initial 
potentials, $\{ f^{XB} \}_B$, into those which approximate the final 
potentials, $\{ f^{AY} \}_A$.
This transition is accomplished in $\log_2 N$ stages, each of which expends at 
most $O(r^2 N^d)$ work in order to map the expansion weights for the 
$N^d$ interactions between members of a target domain partition, $P_X$, and a 
source domain partition, $P_Y$, into weights which approximate the $N^d$ 
interactions between members of a {\em refined} target domain partition, 
$P'_X$, and a {\em coarsened} source domain partition, $P'_Y$.
In particular, $N$ is typically chosen to be a power of two such that, after 
$\ell$ stages, each dimension of the target domain is partitioned into 
$2^\ell$ equal-sized intervals, and each dimension of the source domain is 
like-wise partitioned into $N/2^\ell$ intervals.
This process is depicted for a simple one-dimensional problem, with $N=8$, 
in Fig.~\ref{fig:1D}, and, from now on, we will use the notation 
$\mathcal{T}_X(\ell)$ and $\mathcal{T}_Y(\ell)$ to refer to the sets of 
subdomains produced by the partitions for stage $\ell$, where the symbol 
$\mathcal{T}$ hints at the fact that these are actually trees.
Note that the root of $\mathcal{T}_X$, $\{X\}$, is at stage 0, whereas the 
root of $\mathcal{T}_Y$, $\{Y\}$, is at stage $\log_2 N$.

\subsection{Merging and splitting}
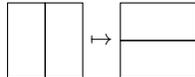
\begin{figure}
\begin{center}
\begin{tikzpicture}[scale=1.0]
\tikzstyle{every node}=[font=\footnotesize]

\draw[black] (0.0,0) rectangle (0.5,1);
\draw[black] (0.5,0) rectangle (1.0,1);

\draw (1.25,0.5) node {$\mapsto$};

\draw[black] (1.5,0.0) rectangle (2.5,0.5);
\draw[black] (1.5,0.5) rectangle (2.5,1.0);

\end{tikzpicture}
\end{center}
\caption{The fundamental operation in the butterfly algorithm: translating from
two potentials with neighboring source domains and equal target domains into 
the corresponding potentials with neighboring target domains but the same 
combined source domain.}
\label{fig:1D-merge}
\end{figure}

A cursory inspection of Fig.~\ref{fig:1D} reveals that each step of a 1D
butterfly algorithm consists of many instances of transforming two potentials 
supported over neighboring source boxes and the same target box, into two 
potentials over the union of the two source boxes but only neighboring halves of
the original target box (see Fig.~\ref{fig:1D-merge}).
The generalization from one to $d$ dimensions is obvious: $2^d$ neighboring 
source boxes are merged and the shared target box is split into $2^d$ 
subboxes.

Suppose that we are given a pair of target boxes, say $A$ and $B$, and we 
define $\{B_j\}_j$ as the set of $2^d$ subboxes of $B$ resulting
from cutting each dimension of the box $B$ into two equal halves.
We shall soon see that, if each pair $(A,B_j)$ satisfies the kernel's 
approximate low-rank criterion, then it is possible to compute a linear 
transformation which, to some predefined accuracy, maps any set of $2^d r$ 
weights representing potentials $\{f^{A B_j}\}_j$ into the 
$2^d r$ weights for the corresponding potentials after the merge-and-split 
procedure, say $\{f^{A_j B}\}_j$, where the $2^d$ subboxes 
$\{A_j\}_j$ of $A$ are defined analogously to those of $B$.

In the case where the source and targets are quasi-uniformly distributed, 
$N^d/2^d$ such linear transformations need to be applied in each of the 
$\log_2 N$ stages of the algorithm, and so, if the corresponding matrices for
these linear transformations have all been precomputed, the per-process cost of
the butterfly algorithm is at most $O(r^2 N^d \log N)$ 
operations~\cite{ONeil-special}.

\subsection{Equivalent sources}
The original approach to the butterfly 
algorithm~\cite{MichielssenBoag-multilevel,ONeil-special} has an elegant 
physical interpretation and provides a straight-forward construction of the 
linear operators which translate weights from one stage to the next.
The primary tool is a (strong) rank-revealing QR (RRQR) 
factorization~\cite{GuEisenstat-RRQR,ChandraIpsen-RRQR}, 
which yields an accurate 
approximation of a numerically low-rank matrix in terms of linear combinations
of a few of its columns.
We will soon see how to manipulate an RRQR into an 
{\em interpolative decomposition} 
(ID)~\cite{ONeil-special,Martinsson-ID},
\begin{equation}
  \mathsf{K} \approx \mathsf{\hat K \hat Z},
\end{equation}
where $\mathsf{\hat K}$ is a submatrix of $r \ll \min(m,n)$ columns of 
the $m \times n$ matrix $\mathsf K$.
While randomized algorithms for building IDs may often be more 
efficient~\cite{Woolfe-randomized}, we will discuss deterministic 
RRQR-based constructions for the sake of simplicity.


Suppose that we have already computed a (truncated) RRQR decomposition
\begin{equation} 
  \mathsf{K \Pi} \mathsf{\approx Q} 
  \begin{pmatrix} \mathsf{R_L} & \mathsf{R_R} \end{pmatrix}, 
\end{equation}
where $\mathsf{\Pi}$ is a permutation matrix, $\mathsf{Q}$ consists of $r$ 
mutually orthonormal columns, and $\mathsf{R_L}$ is an invertible $r \times r$ 
upper-triangular matrix (otherwise a smaller RRQR could be formed).
Then 
$\mathsf{K \Pi} \approx \mathsf{\hat K} 
 \begin{pmatrix} \mathsf{I} & \mathsf{R_L^{-1} R_R} \end{pmatrix}$, 
where $\mathsf{\hat K} \equiv \mathsf{Q R_L}$ is the subset of columns of $\mathsf K$ 
selected during the pivoted QR factorization (the first $r$ columns of 
$\mathsf{K \Pi}$), and 
\begin{equation}
  \mathsf{\hat Z} \equiv \begin{pmatrix} \mathsf{I} & \mathsf{R_L^{-1} R_R} \end{pmatrix} \mathsf{\Pi^H}
\end{equation}
is an $r \times n$ interpolation matrix such that 
$\mathsf K \approx \mathsf{\hat K \hat Z}$, which completes our interpolative 
decomposition of $\mathsf K$.
We note that, although $\mathsf{\hat Z}$ is only guaranteed to be computed 
stably from a {\em strong} RRQR factorization, Businger-Golub 
pivoting~\cite{BusingerGolub} is typically used in practice in combination 
with applying the numerical pseudoinverse of $R_L$ rather than 
$R_L^{-1}$~\cite{Martinsson-ID}.

The algorithms of \cite{MichielssenBoag-multilevel,ONeil-special} 
exploit the fact that the interpolation matrix, $\mathsf{\hat Z}$, can be 
applied to a dense vector of $n$ ``sources'' in order to produce an 
approximately equivalent set of $r$ sources, in the sense that, for any vector 
$\mathsf{g} \in \mathbb{C}^n$, 
\begin{equation}
  \mathsf{K g \approx \hat K (\hat Z g) = \hat K \hat g}. 
\end{equation}
Or, equivalently, $\mathsf{K g} \approx \mathsf{K \hat{g}_E}$, 
where $\mathsf{\hat{g}_E}$ is the appropriate extension by zero of 
$\mathsf{\hat g} \in \mathbb{C}^r$ into $\mathbb{C}^n$.
$\mathsf{\hat Z}$ therefore provides a fast mechanism for producing an 
approximately equivalent \emph{sparse} source vector, $\mathsf{\hat{g}_E}$, 
given any (potentially dense) source vector, $\mathsf g$.

These sparse (approximately) equivalent sources can then be used as the 
expansion weights resulting from the low-rank approximation 
\[
  \left| 
  \mathsf{K(i,j)} - \sum_{t=0}^{r-1} \mathsf{K(i,j_t) \hat z_t(j)} 
  \right| \le \epsilon\, s(r,n),\;\;\;\forall\,i,j,
\]
namely,
\[
  \left\| \mathsf{K g} - \mathsf{\hat{K} \hat{g}} \right\|_\infty \le 
  \epsilon\, s(r,n) \| \mathsf{g} \|_1,
\]
where $\mathsf{K(:,j_t)}$ is the $t$'th column of $\mathsf{\hat K}$, 
$\mathsf{\hat z_t}$ is the $t$'th row of $\mathsf{\hat Z}$, 
$\mathsf{\hat{g}} \equiv \mathsf{\hat{Z} g}$, and $s(r,n)$, 
which is bounded by a low-degree polynomial in $r$ and $n$, is an artifact of 
RRQR factorizations yielding suboptimal low-rank 
decompositions~\cite{GuEisenstat-RRQR}.

If the matrix $\mathsf{K}$ and vector $\mathsf{g}$ were constructed such that 
$\mathsf{K(i,j)}=K(x_i,y_j)$ and $\mathsf{g}(j)=g(y_j)$, for some set of 
source points $\{y_j\}_j \subset B$ and target points $\{x_i\}_i \subset A$, 
then the previous equation becomes
\[
  \left| f^{AB}(x_i) - \sum_{t=0}^{r-1} K(x_i,y_{j_t}) \hat g(y_{j_t}) \right| 
  \le \epsilon\, s(r,n) \| g_B \|_1,\;\;\; \forall\, i,
\]
where we have emphasized the interpretation of the $t$'th entry of 
the equivalent source vector $\mathsf{\hat{g}}$ as a discrete source located 
at the point $y_{j_t}$, i.e., $\hat g(y_{j_t})$.
We will now review how repeated applications of IDs can yield
interpolation matrices which take $2^d r$ sources from $2^d$ neighboring source
boxes and produce (approximately) equivalent sets of $r$ sources valid over
smaller target domains. 

\subsection{Translating equivalent sources}
We will begin by considering the one-dimensional case, as it lends itself to a 
matrix-centric discussion:
Let $B_0$ and $B_1$ be two neighboring source intervals of the same size, let 
$A$ be a target interval of sufficiently small width, and let 
$\mathsf{\hat K_{A B_0}}$ and $\mathsf{\hat K_{A B_1}}$ be subsets of 
columns from $\mathsf{K|_{A B_0}}$ and $\mathsf{K|_{A B_1}}$ generated from 
their interpolative decompositions with interpolation matrices 
$\mathsf{\hat Z_{A B_0}}$ and $\mathsf{\hat Z_{A B_1}}$.
Then, for any source vector $\mathsf g$, 
the vector of potentials over target box $A$ generated by the sources in box 
$B_n$ can be cheaply approximated as
\[
  \mathsf{f_{A B_n}} \approx \mathsf{\hat K_{A B_n}}\, \mathsf{\hat g_{A B_n}} 
  \left(=\mathsf{\hat K_{A B_n}\, \hat Z_{A B_n}\, g|_{B_n}}\right),\; 
  n=0,1.
\]
If we then define $B = B_0 \cup B_1$ and halve the interval $A$ into 
$A_0$ and $A_1$, then the products of the widths of each $A_m$ with 
the entire box $B$ is equal to that of $A$ with each $B_n$, and so, due to the 
main assumption of the butterfly algorithm, $\mathsf{K|_{A_m B}}$ must also be 
numerically low-rank.

If we then split each $\mathsf{\hat K_{A B_n}}$ into the two submatrices
\[
  \mathsf{\hat K_{A B_n}} \rightarrow 
   \begin{bmatrix} \mathsf{\hat K_{A_0 B_n}} \\ 
                   \mathsf{\hat K_{A_1 B_n}} \end{bmatrix},
\]
we can write
\[
  \begin{pmatrix}
  \mathsf{f_{A_0 B}} \\
  \mathsf{f_{A_1 B}} 
  \end{pmatrix}
  \approx 
  \begin{bmatrix} 
  \mathsf{\hat K_{A_0 B_0}} & \mathsf{\hat K_{A_0 B_1}} \\
  \mathsf{\hat K_{A_1 B_0}} & \mathsf{\hat K_{A_1 B_1}} 
  \end{bmatrix}
  \begin{pmatrix}
  \mathsf{\hat g_{A B_0}} \\
  \mathsf{\hat g_{A B_1}}
  \end{pmatrix} 
\]
and recognize that the two submatrices 
$\mathsf{\left[ \hat K_{A_m B_0}, \hat K_{A_m B_1} \right]}$, $m=0,1$, 
consist of subsets of columns of $\mathsf{K|_{A_m B}}$, which implies that 
they must also have low-rank interpolative decompositions, 
say $\mathsf{\hat K_{A_m B} \hat Z_{A_m B}}$.
Thus,
\[
  \begin{pmatrix}
  \mathsf{f_{A_0 B}} \\
  \mathsf{f_{A_1 B}} 
  \end{pmatrix}
  \approx 
  \begin{bmatrix} 
  \mathsf{\hat K_{A_0 B} \hat Z_{A_0 B}} \\
  \mathsf{\hat K_{A_1 B} \hat Z_{A_1 B}} 
  \end{bmatrix}
  \begin{pmatrix}
  \mathsf{\hat g_{A B_0}} \\
  \mathsf{\hat g_{A B_1}}
  \end{pmatrix} 
  =
  \begin{pmatrix} 
  \mathsf{\hat K_{A_0 B} \hat g_{A_0 B}} \\
  \mathsf{\hat K_{A_1 B} \hat g_{A_1 B}}
  \end{pmatrix},
\]
where we have defined the new equivalent sources, $\mathsf{\hat g_{A_m B}}$, as
\[
  \mathsf{\hat g_{A_m B} \equiv \hat Z_{A_m B}} 
  \begin{pmatrix} 
    \mathsf{\hat g_{A B_0}} \\ \mathsf{\hat g_{A B_1}}
  \end{pmatrix}.
\]

The generalization to $d$-dimensions should again be clear: $2^d$ 
sets of IDs should be stacked together, partitioned, and recompressed in 
order to form $2^d$ interpolation matrices of size $r \times 2^d r$, 
say $\{\mathsf{\hat Z_{A_m B}}\}_m$. 
These interpolation matrices can then be used to translate $2^d$ sets of 
equivalent sources from one stage to the next with at most $O(r^2)$ work.
Recall that each of the $\log_2 N$ stage of the butterfly algorithm requires 
$O(N^d)$ such translations, and so, \emph{if all necessary IDs have been 
precomputed}, the equivalent source approach yields an $O(r^2 N^d \log N)$ 
butterfly algorithm. 
There is, of course, an analogous approach based on row-space interpolation, 
which can be interpreted as constructing a small set of representative 
potentials which may then be cheaply interpolated to evaluate the potential 
over the entire target box.

Yet another approach would be to replace row/column-space interpolation with
the low-rank approximation implied by a Singular Value Decomposition (SVD) and 
to construct the $2^d r \times 2^d r$ linear weight transformation matrices
based upon the low-rank approximations used at successive levels.
Since the low-rank approximations produced by SVDs are generally much tighter
than those of rank-revealing factorizations, such an approach would potentially
allow for lower-rank approximations to result in the same overall 
accuracy.\footnote{A small amount of structure in the translation operators is forfeited 
when switching to an SVD-based approach: the left-most $r \times r$ subblock of the 
(permuted) $r \times 2^d r$ translation matrix changes from the identity to an arbitrary 
(dense) matrix.}

From now on, we will use the high-level notation that 
$\mathsf{T_{A_m B}^{A\, B_n}}$ is the {\em translation operator} which maps 
a weight vector from a low-rank decomposition over $A \times B_n$, 
$\mathsf{w_{A B_n}}$, into that of a low-rank decomposition over 
$A_m \times B$, $\mathsf{w_{A_m B}}$.
Clearly each merge and split operation involves a $2^d \times 2^d$ block matrix
of such translation operators.
Please see Algorithm~\ref{alg:serial} for a demonstration of the sequential 
algorithm from the point of view of translation operators, where the low-rank
approximation of each block $\mathsf{K|_{AB}}$ is denoted by 
$\mathsf{U_{AB} V_{AB}}$, and we recall that the nodes of the trees active at 
the beginning of the $\ell$'th stage of the algorithm are denoted by 
$\mathcal{T}_X(\ell)$ and $\mathcal{T}_Y(\ell)$.

\begin{algorithm}
\DontPrintSemicolon
$\mathcal{A} := \mathcal{T}_X(0),\; \mathcal{B} := \mathcal{T}_Y(0)$\;
$\left(\mathcal{A} = \{X\},\; 
  \cup_{B\in \mathcal{B}} B = Y,\; \text{card}(\mathcal{B})=N^d\right)$\;
\tcp{Initialize weights}
\ForEach{$B \in \mathcal{B}$}{
  $\mathsf{w_{XB} := V_{XB}\, g_B}$\;
}
\tcp{Translate weights}
\For{$\ell=0,...,\log_2 N - 1$}{
  $\tilde{\mathcal{A}} := \text{children}(\mathcal{A}),\;
   \tilde{\mathcal{B}} := \text{parents}(\mathcal{B})$\;
  \ForEach{$(\tilde A,\tilde B) \in 
            \tilde{\mathcal{A}} \times \tilde{\mathcal{B}}$}{
    $\mathsf{w_{\tilde A \tilde B} := 0}$\;
  }
  \ForEach{$(A,B) \in \mathcal{A} \times \mathcal{B}$}{
    $\{A_c\}_{c=0}^{2^d-1} := \text{children}(A),\; B_p := \text{parent}(B)$\;
    \ForEach{$c=0,...,2^d-1$}{
      $\mathsf{w_{A_c B_p}\, +\!\!= T_{A_c B_p}^{A B} w_{AB}}$\;
    }
  }
  $\mathcal{A} := \tilde{\mathcal{A}},\; \mathcal{B} := \tilde{\mathcal{B}}$\;
}
$\left(\text{card}(\mathcal{A})=N^d,\; 
  \cup_{A\in \mathcal{A}} A = X,\; \mathcal{B} = \{Y\}\right)$\;
\tcp{Final evaluations}
\ForEach{$A \in \mathcal{A}$}{
  $\mathsf{f_{AY} := U_{AY}\, w_{AY}}$\;
}
\caption{Sequential butterfly algorithm over a $d$-dimensional domain with
$N^d$ source and target points.}
\label{alg:serial}
\end{algorithm}

\subsection{Avoiding quadratic precomputation}
The obvious drawback to ID and SVD-based approaches is that the precomputation 
of the $O(N^d \log N)$ necessary low-rank approximations requires at least 
$O(N^{2d})$ work with any black-box algorithm, as the first stage of the 
butterfly algorithm involves $O(N^d)$ matrices of height $N^d$.
If we assume additional features of the underlying kernel, we may 
accelerate these precomputations~\cite{Tygert-spherical,Seljebotn-wavemoth} or, 
in some cases, essentially avoid them altogether~\cite{CDY-butterfly,Ying-SFT}.

We focus on the latter case, where the kernel is assumed to be of the form
\begin{equation}\label{phase}
  K(x,y) = e^{i \Phi(x,y)},
\end{equation}
where $\Phi : X \times Y \rightarrow \mathbb{R}$ is 
a sufficiently smooth\footnote{Formally, (Q,R)-analytic~\cite{CDY-butterfly,Demanet-SAR}} \emph{phase function}.
Due to the assumed smoothness of $\Phi$, it was shown in 
\cite{CDY-butterfly,Demanet-SAR,KunisMelzer-stableButterfly} that the 
precomputation of IDs can be replaced with analytical interpolation of the 
row or column space of the numerically low-rank submatrices using tensor-product
Chebyshev grids.
In particular, in the first half of the algorithm, while the source boxes are
small, interpolation is performed within the column space, and in the middle 
of the algorithm, when the target boxes become as small as the source boxes, 
the algorithm switches to analytical row-space interpolation.

An added benefit of the tensor-product interpolation is that, if a basis of 
dimension $q$ is used in each direction, so that the rank of each 
approximation is $r=q^d$, the weight translation cost can be reduced from 
$O(r^2)$ to $O(q^{d+1})=O(r^{1+1/d})$. However, we note that the cost of the 
{\em middle-switch} from column-space to row-space interpolation, in general,
requires $O(r^2)$ work for each set  of weights. But this cost can also be 
reduced to $O(r^{1+1/d})$ when the phase function, $\Phi(x,y)$, also has a 
tensor-product structure, e.g., $x \cdot y$.
Lastly, performing analytical interpolation allows one to choose the precise 
locations of the target points {\em after} forming the final expansion weights,
$\{\mathsf{w_{AY}}\}_A$, and so the result is best viewed as a potential field 
which only requires $O(r)$ work to evaluate at any given point in the 
continuous target domain.

\section{Parallelization}
We now present a parallelization of butterfly algorithms which is high-level 
enough to handle both 
general-purpose~\cite{MichielssenBoag-multilevel,ONeil-special} and 
analytical~\cite{CDY-butterfly} low-rank interpolation.
We will proceed by first justifying our communication cost model, then step 
through the simplest parallel case, where each box is responsible for a single 
interaction at a time, and then demonstrate how the communication requirements 
change when less processes are used.
We will not discuss the precomputation phase in detail, as the factorizations
computed at each level need to be exchanged between processes in the same 
manner as the weights, but the factorizations themselves are relatively much 
more expensive and can each be run sequentially (unless $p>N^d$). 
Readers interested in efficient parallel IDs should consult the 
communication-avoiding RRQR factorization of \cite{Demmel-CARRQR}.

\subsection{Communication cost model}
All of our analysis makes use of a 
commonly-used~\cite{FosterWorley-STM,TRG-coll,Chan-collective,Ballard-MinComm,Demmel-CARRQR} 
communication cost model that is as useful as it is simple: each process is 
assumed to only be able to simultaneously send and receive a single message at 
a time, and, when the message consists of $n$ units of data 
(e.g., double-precision floating-point numbers), 
the time to transmit such a message between any two processes is 
$\alpha + \beta n$~\cite{Hockney,InterCom}.
The $\alpha$ term represents the time required to send an arbitrarily small 
message and is commonly referred to as the message {\em latency}, whereas
$1/\beta$ represents the number of units of data which can be transmitted 
per unit of time once the message has been initiated.

There also exist more sophisticated communication models, such as
{\em LogP}~\cite{LogP} and its extension, {\em LogGP}~\cite{LogGP}, but the 
essential differences are that the former separates the local software 
overhead (the `o' in `LogP') from the network latency and the latter 
compensates for very large messages potentially having a different 
transmission mechanism.
An arguably more important detail left out of the $\alpha + \beta n$ model
is the issue of {\em network conflicts}~\cite{TRG-coll,Chan-collective}, 
that is, when multiple messages compete for the bandwidth available on a 
single communication link.
We will ignore network conflicts since they greatly complicate our analysis
and require the specialization of our cost model to a particular network 
topology.

\subsection{High-level approach with $N^d$ processes}
The proposed parallelization is easiest to discuss for cases where the 
number of processes is equal to $N^d$, the number of pairwise source and 
target box interactions represented at each level of the butterfly algorithm.
Recall that each of these interactions is represented with $r$ expansion 
weights, where $r$ is a small number which should depend polylogarithmically
on the desired accuracy.
We will now present a scheme which assigns each process one set of expansion
weights at each stage and has a runtime of $O((r^2 + \beta r + \alpha) \log N)$.


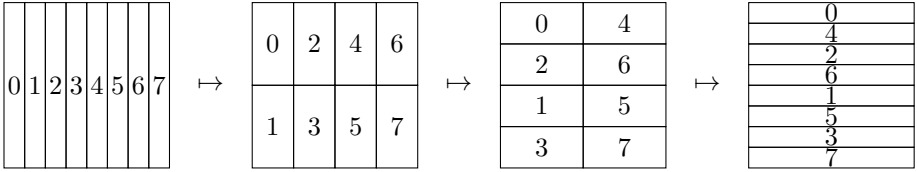
\begin{figure}
\begin{center}
\begin{tikzpicture}[scale=2.2]

\foreach \j in {0,...,7}{
  \draw[black] (0.125*\j,0) rectangle (0.125*\j+0.125,1);
  \draw (0.125*\j+0.0625,0.5) node {$\j$};
}

\draw (1.25,0.5) node {$\mapsto$};

\draw[black] (1.5+0.25*0,0.5*0) rectangle (1.5+0.25*0+0.25,0.5*0+0.5);
\draw[black] (1.5+0.25*1,0.5*0) rectangle (1.5+0.25*1+0.25,0.5*0+0.5);
\draw[black] (1.5+0.25*2,0.5*0) rectangle (1.5+0.25*2+0.25,0.5*0+0.5);
\draw[black] (1.5+0.25*3,0.5*0) rectangle (1.5+0.25*3+0.25,0.5*0+0.5);
\draw[black] (1.5+0.25*0,0.5*1) rectangle (1.5+0.25*0+0.25,0.5*1+0.5);
\draw[black] (1.5+0.25*1,0.5*1) rectangle (1.5+0.25*1+0.25,0.5*1+0.5);
\draw[black] (1.5+0.25*2,0.5*1) rectangle (1.5+0.25*2+0.25,0.5*1+0.5);
\draw[black] (1.5+0.25*3,0.5*1) rectangle (1.5+0.25*3+0.25,0.5*1+0.5);
\draw (1.625+0.25*0,0.25+0.5*0) node {$1$};
\draw (1.625+0.25*1,0.25+0.5*0) node {$3$};
\draw (1.625+0.25*2,0.25+0.5*0) node {$5$};
\draw (1.625+0.25*3,0.25+0.5*0) node {$7$};
\draw (1.625+0.25*0,0.25+0.5*1) node {$0$};
\draw (1.625+0.25*1,0.25+0.5*1) node {$2$};
\draw (1.625+0.25*2,0.25+0.5*1) node {$4$};
\draw (1.625+0.25*3,0.25+0.5*1) node {$6$};

\draw (2.75,0.5) node {$\mapsto$};

\draw[black] (3+0.5*0,0.25*0) rectangle (3+0.5*0+0.5,0.25*0+0.25);
\draw[black] (3+0.5*1,0.25*0) rectangle (3+0.5*1+0.5,0.25*0+0.25);
\draw (3.25+0.5*0,0.125+0.25*0) node {$3$};
\draw (3.25+0.5*1,0.125+0.25*0) node {$7$};

\draw[black] (3+0.5*0,0.25*1) rectangle (3+0.5*0+0.5,0.25*1+0.25);
\draw[black] (3+0.5*1,0.25*1) rectangle (3+0.5*1+0.5,0.25*1+0.25);
\draw (3.25+0.5*0,0.125+0.25*1) node {$1$};
\draw (3.25+0.5*1,0.125+0.25*1) node {$5$};

\draw[black] (3+0.5*0,0.25*2) rectangle (3+0.5*0+0.5,0.25*2+0.25);
\draw[black] (3+0.5*1,0.25*2) rectangle (3+0.5*1+0.5,0.25*2+0.25);
\draw (3.25+0.5*0,0.125+0.25*2) node {$2$};
\draw (3.25+0.5*1,0.125+0.25*2) node {$6$};

\draw[black] (3+0.5*0,0.25*3) rectangle (3+0.5*0+0.5,0.25*3+0.25);
\draw[black] (3+0.5*1,0.25*3) rectangle (3+0.5*1+0.5,0.25*3+0.25);
\draw (3.25+0.5*0,0.125+0.25*3) node {$0$};
\draw (3.25+0.5*1,0.125+0.25*3) node {$4$};

\draw (4.25,0.5) node {$\mapsto$};

\draw[black] (4.5+0,0.125*0) rectangle (4.5+1,0.125*0+0.125);
\draw[black] (4.5+0,0.125*1) rectangle (4.5+1,0.125*1+0.125);
\draw[black] (4.5+0,0.125*2) rectangle (4.5+1,0.125*2+0.125);
\draw[black] (4.5+0,0.125*3) rectangle (4.5+1,0.125*3+0.125);
\draw[black] (4.5+0,0.125*4) rectangle (4.5+1,0.125*4+0.125);
\draw[black] (4.5+0,0.125*5) rectangle (4.5+1,0.125*5+0.125);
\draw[black] (4.5+0,0.125*6) rectangle (4.5+1,0.125*6+0.125);
\draw[black] (4.5+0,0.125*7) rectangle (4.5+1,0.125*7+0.125);
\draw (4.5+0.5,0.0625+0.125*0) node {$7$};
\draw (4.5+0.5,0.0625+0.125*1) node {$3$};
\draw (4.5+0.5,0.0625+0.125*2) node {$5$};
\draw (4.5+0.5,0.0625+0.125*3) node {$1$};
\draw (4.5+0.5,0.0625+0.125*4) node {$6$};
\draw (4.5+0.5,0.0625+0.125*5) node {$2$};
\draw (4.5+0.5,0.0625+0.125*6) node {$4$};
\draw (4.5+0.5,0.0625+0.125*7) node {$0$};

\end{tikzpicture}
\end{center}
\caption{The data distribution throughout the execution of a 1D butterfly 
algorithm with $p=N=8$. Notice that the distribution of the target domain upon 
output is the same as the input distribution of the source domain, but with 
the binary process ranks reversed.}
\label{fig:1D-par}
\end{figure}

Consider the data distribution scheme shown in Fig.~\ref{fig:1D-par} for a 
1D butterfly algorithm with both the problem size, $N$, and number of 
processes, $p$, set to eight.
In the beginning of the algorithm (the left-most panel), each process is 
assigned one source box and need only expend $O(r)$ work in order to 
initialize its weights, and, at the end of the algorithm (the right-most panel),
each process can perform $O(r)$ flops in order to evaluate the potential over 
its target box.
Since each stage of the butterfly algorithm involves linearly transforming 
$2^d$ sets of weights from one level to the next, e.g., via 
\[
  \begin{pmatrix} 
    \mathsf{w_{A_0 B}} \\ \mathsf{w_{A_1 B}} 
  \end{pmatrix} = 
  \begin{bmatrix}[1.5]
    \mathsf{T_{A_0 B}^{A\, B_0}} & \mathsf{T_{A_0 B}^{A\, B_1}} \\
    \mathsf{T_{A_1 B}^{A\, B_0}} & \mathsf{T_{A_1 B}^{A\, B_1}} 
  \end{bmatrix}
  \begin{pmatrix}
    \mathsf{w_{A B_0}} \\
    \mathsf{w_{A B_1}}
  \end{pmatrix},
\]
pairs of processes will need to coordinate in order to perform 
parallel matrix-vector multiplications of size $2 r \times 2 r$.
Because each process initially owns only half of the vector that must be 
linearly transformed, 
it is natural to locally compute half of the linear transformation, for example,
$\mathsf{T_{A_0 B}^{A\, B_0}\, w_{A B_0}}$ and 
$\mathsf{T_{A_1 B}^{A B_0}\, w_{A B_0}}$, and to combine the results with those 
computed by the partner process.
Furthermore, the output weights of such a linear transformation should
also be distributed, and so only $r$ entries of data need be exchanged between 
the two processes, for a cost of $\alpha + \beta r$.
Since the cost of the local transformation is at most $O(r^2)$, and only 
$O(r)$ work is required to sum the received data, the cost of the each stage is 
at most $O(r^2 + \beta r + \alpha)$.

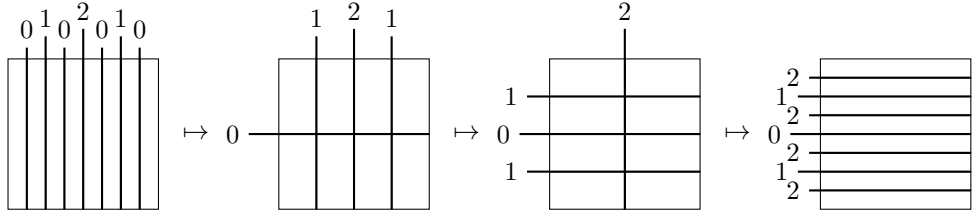
\begin{figure}
\begin{center}
\begin{tikzpicture}[scale=2]

\foreach \j in {0,...,7}{
  \draw[black] (0.125*\j,0) rectangle (0.125*\j+0.125,1);
}
\draw[black,thick] (0.125*4,0) -- (0.125*4,1.2);
\draw (0.125*4,1.2) node[above] {$2$};

\draw[black,thick] (0.125*2,0) -- (0.125*2,1.15);
\draw[black,thick] (0.125*6,0) -- (0.125*6,1.15);
\draw (0.125*2,1.15) node[above] {$1$};
\draw (0.125*6,1.15) node[above] {$1$};

\draw[black,thick] (0.125*1,0) -- (0.125*1,1.075);
\draw[black,thick] (0.125*3,0) -- (0.125*3,1.075);
\draw[black,thick] (0.125*5,0) -- (0.125*5,1.075);
\draw[black,thick] (0.125*7,0) -- (0.125*7,1.075);
\draw (0.125*1,1.075) node[above] {$0$};
\draw (0.125*3,1.075) node[above] {$0$};
\draw (0.125*5,1.075) node[above] {$0$};
\draw (0.125*7,1.075) node[above] {$0$};

\draw (1.25,0.5) node {$\mapsto$};

\draw[black] (1.8+0.25*0,0.5*0) rectangle (1.8+0.25*0+0.25,0.5*0+0.5);
\draw[black] (1.8+0.25*1,0.5*0) rectangle (1.8+0.25*1+0.25,0.5*0+0.5);
\draw[black] (1.8+0.25*2,0.5*0) rectangle (1.8+0.25*2+0.25,0.5*0+0.5);
\draw[black] (1.8+0.25*3,0.5*0) rectangle (1.8+0.25*3+0.25,0.5*0+0.5);
\draw[black] (1.8+0.25*0,0.5*1) rectangle (1.8+0.25*0+0.25,0.5*1+0.5);
\draw[black] (1.8+0.25*1,0.5*1) rectangle (1.8+0.25*1+0.25,0.5*1+0.5);
\draw[black] (1.8+0.25*2,0.5*1) rectangle (1.8+0.25*2+0.25,0.5*1+0.5);
\draw[black] (1.8+0.25*3,0.5*1) rectangle (1.8+0.25*3+0.25,0.5*1+0.5);

\draw[black,thick] (1.8+0.125*4,0) -- (1.8+0.125*4,1.2);
\draw (1.8+0.125*4,1.2) node[above] {$2$};

\draw[black,thick] (1.8+0.125*2,0) -- (1.8+0.125*2,1.15);
\draw[black,thick] (1.8+0.125*6,0) -- (1.8+0.125*6,1.15);
\draw (1.8+0.125*2,1.15) node[above] {$1$};
\draw (1.8+0.125*6,1.15) node[above] {$1$};

\draw[thick] (1.6,0.5) -- (2.8,0.5);
\draw (1.6,0.5) node[left] {$0$};

\draw (3.05,0.5) node {$\mapsto$};

\draw[black] (3.6+0.5*0,0.25*0) rectangle (3.6+0.5*0+0.5,0.25*0+0.25);
\draw[black] (3.6+0.5*1,0.25*0) rectangle (3.6+0.5*1+0.5,0.25*0+0.25);

\draw[black] (3.6+0.5*0,0.25*1) rectangle (3.6+0.5*0+0.5,0.25*1+0.25);
\draw[black] (3.6+0.5*1,0.25*1) rectangle (3.6+0.5*1+0.5,0.25*1+0.25);

\draw[black] (3.6+0.5*0,0.25*2) rectangle (3.6+0.5*0+0.5,0.25*2+0.25);
\draw[black] (3.6+0.5*1,0.25*2) rectangle (3.6+0.5*1+0.5,0.25*2+0.25);

\draw[black] (3.6+0.5*0,0.25*3) rectangle (3.6+0.5*0+0.5,0.25*3+0.25);
\draw[black] (3.6+0.5*1,0.25*3) rectangle (3.6+0.5*1+0.5,0.25*3+0.25);

\draw[black,thick] (3.6+0.125*4,0) -- (3.6+0.125*4,1.2);
\draw (3.6+0.125*4,1.2) node[above] {$2$};

\draw[thick] (3.45,0.25) -- (4.6,0.25);
\draw[thick] (3.45,0.75) -- (4.6,0.75);
\draw (3.45,0.25) node[left] {$1$};
\draw (3.45,0.75) node[left] {$1$};

\draw[thick] (3.4,0.5) -- (4.6,0.5);
\draw (3.4,0.5) node[left] {$0$};

\draw (4.85,0.5) node {$\mapsto$};

\draw[black] (5.4+0,0.125*0) rectangle (5.4+1,0.125*0+0.125);
\draw[black] (5.4+0,0.125*1) rectangle (5.4+1,0.125*1+0.125);
\draw[black] (5.4+0,0.125*2) rectangle (5.4+1,0.125*2+0.125);
\draw[black] (5.4+0,0.125*3) rectangle (5.4+1,0.125*3+0.125);
\draw[black] (5.4+0,0.125*4) rectangle (5.4+1,0.125*4+0.125);
\draw[black] (5.4+0,0.125*5) rectangle (5.4+1,0.125*5+0.125);
\draw[black] (5.4+0,0.125*6) rectangle (5.4+1,0.125*6+0.125);
\draw[black] (5.4+0,0.125*7) rectangle (5.4+1,0.125*7+0.125);

\draw[black,thick] (5.325,0.125) -- (6.4,0.125);
\draw[black,thick] (5.325,0.375) -- (6.4,0.375);
\draw[black,thick] (5.325,0.625) -- (6.4,0.625);
\draw[black,thick] (5.325,0.875) -- (6.4,0.875);
\draw (5.325,0.125) node[left] {$2$};
\draw (5.325,0.375) node[left] {$2$};
\draw (5.325,0.625) node[left] {$2$};
\draw (5.325,0.875) node[left] {$2$};

\draw[thick] (5.25,0.25) -- (6.4,0.25);
\draw[thick] (5.25,0.75) -- (6.4,0.75);
\draw (5.25,0.25) node[left] {$1$};
\draw (5.25,0.75) node[left] {$1$};

\draw[thick] (5.2,0.5) -- (6.4,0.5);
\draw (5.2,0.5) node[left] {$0$};

\end{tikzpicture}
\end{center}
\caption{An alternative view of the data distributions used during a 1D 
parallel butterfly algorithm with $p=N=8$ based upon bitwise bisections:
processes with their $j$'th bit set to zero are assigned to the left or 
upper side of the partition, while processes with the $j$'th bit of their 
rank set to one are assigned to the other side.}
\label{fig:1D-par-bin}
\end{figure}

In the first, second, and third stages of the algorithm, process 0 would 
respectively pair with processes 1, 2, and 4, and its combined communication and
computation cost would be $O((r^2 + \beta r + \alpha) \log N)$. 
Every other process need only perform the same amount of work and communication,
and it can be seen from the right-most panel of Fig.~\ref{fig:1D-par} that,
upon completion, each process will hold an approximation for the potential 
generated over a single target box resulting from the entire set of sources.
In fact, the final distribution of the target domain can be seen to be the same
as the initial distribution of the source domain, but with the bits of the 
owning processes reversed (see Fig.~\ref{fig:1D-par-bin}).

\begin{figure}
\begin{center}
\begin{tikzpicture}[scale=1.8]

\foreach \j in {0,...,3}{
  \draw[black] (0.25*\j,0) rectangle (0.25*\j+0.25,1);
}
\draw[black,thick] (0.125*4,0) -- (0.125*4,1.2);
\draw (0.125*4,1.2) node[above] {$5$};

\draw[black,thick] (0.125*2,0) -- (0.125*2,1.15);
\draw[black,thick] (0.125*6,0) -- (0.125*6,1.15);
\draw (0.125*2,1.15) node[above] {$3$};
\draw (0.125*6,1.15) node[above] {$3$};

\draw[black,thick] (0.125*1,0) -- (0.125*1,1.075);
\draw[black,thick] (0.125*3,0) -- (0.125*3,1.075);
\draw[black,thick] (0.125*5,0) -- (0.125*5,1.075);
\draw[black,thick] (0.125*7,0) -- (0.125*7,1.075);
\draw (0.125*1,1.075) node[above] {$1$};
\draw (0.125*3,1.075) node[above] {$1$};
\draw (0.125*5,1.075) node[above] {$1$};
\draw (0.125*7,1.075) node[above] {$1$};

\draw[black,thick] (-0.2,0.5) -- (1,0.5);
\draw (-0.2,0.5) node[left] {$4$};

\draw[black,thick] (-0.15,0.25) -- (1,0.25);
\draw[black,thick] (-0.15,0.75) -- (1,0.75);
\draw (-0.15,0.25) node[left] {$2$};
\draw (-0.15,0.75) node[left] {$2$};

\draw[black,thick] (-0.075,1*0.125) -- (1,1*0.125);
\draw[black,thick] (-0.075,3*0.125) -- (1,3*0.125);
\draw[black,thick] (-0.075,5*0.125) -- (1,5*0.125);
\draw[black,thick] (-0.075,7*0.125) -- (1,7*0.125);
\draw[black,thick] (-0.075,1*0.125) node[left] {$0$};
\draw[black,thick] (-0.075,3*0.125) node[left] {$0$};
\draw[black,thick] (-0.075,5*0.125) node[left] {$0$};
\draw[black,thick] (-0.075,7*0.125) node[left] {$0$};

\draw[black] (0,1.65) rectangle (1,2.65);

\draw (1.4,1.375) node {\Large $\mapsto$};

\draw[black] (1.8+0.5*0,0.5*0) rectangle (1.8+0.5*1,0.5*0+0.5);
\draw[black] (1.8+0.5*1,0.5*0) rectangle (1.8+0.5*2,0.5*0+0.5);
\draw[black] (1.8+0.5*0,0.5*1) rectangle (1.8+0.5*1,0.5*1+0.5);
\draw[black] (1.8+0.5*1,0.5*1) rectangle (1.8+0.5*2,0.5*1+0.5);

\draw[black,thick] (1.8+0.125*4,0) -- (1.8+0.125*4,1.2);
\draw (1.8+0.125*4,1.2) node[above] {$5$};

\draw[black,thick] (1.8+0.125*2,0) -- (1.8+0.125*2,1.15);
\draw[black,thick] (1.8+0.125*6,0) -- (1.8+0.125*6,1.15);
\draw (1.8+0.125*2,1.15) node[above] {$3$};
\draw (1.8+0.125*6,1.15) node[above] {$3$};

\draw[black,thick] (1.6,0.5) -- (2.8,0.5);
\draw (1.6,0.5) node[left] {$4$};

\draw[black,thick] (1.65,0.25) -- (2.8,0.25);
\draw[black,thick] (1.65,0.75) -- (2.8,0.75);
\draw (1.65,0.25) node[left] {$2$};
\draw (1.65,0.75) node[left] {$2$};

\draw[black] (1.8+0.5*0,1.65+0.5*0) rectangle (1.8+0.5*1,0.5*0+2.15);
\draw[black] (1.8+0.5*1,1.65+0.5*0) rectangle (1.8+0.5*2,0.5*0+2.15);
\draw[black] (1.8+0.5*0,1.65+0.5*1) rectangle (1.8+0.5*1,0.5*1+2.15);
\draw[black] (1.8+0.5*1,1.65+0.5*1) rectangle (1.8+0.5*2,0.5*1+2.15);

\draw[black,thick] (1.8+0.125*4,1.65) -- (1.8+0.125*4,2.85);
\draw (1.8+0.125*4,2.85) node[above] {$0$};

\draw[black,thick] (1.6,2.15) -- (2.8,2.15);
\draw (1.6,2.15) node[left] {$1$};

\draw (3.2,1.375) node {\Large $\mapsto$};

\draw[black] (3.6,0) rectangle (4.6,1);

\draw[black,thick] (3.6+0.125*4,0) -- (3.6+0.125*4,1.2);
\draw (3.6+0.125*4,1.2) node[above] {$5$};

\draw[black,thick] (3.4,0.5) -- (4.6,0.5);
\draw (3.4,0.5) node[left] {$4$};

\foreach \j in {0,...,3}{
  \draw[black] (3.6+0.25*\j,1.65) rectangle (3.6+0.25*\j+0.25,2.65);
}
\draw[black,thick] (3.6+0.125*4,1.65) -- (3.6+0.125*4,2.85);
\draw (3.6+0.125*4,2.85) node[above] {$0$};

\draw[black,thick] (3.6+0.125*2,1.65) -- (3.6+0.125*2,2.75);
\draw[black,thick] (3.6+0.125*6,1.65) -- (3.6+0.125*6,2.75);
\draw (3.6+0.125*2,2.75) node[above] {$2$};
\draw (3.6+0.125*6,2.75) node[above] {$2$};

\draw[black,thick] (3.4,2.15) -- (4.6,2.15);
\draw (3.4,2.15) node[left] {$1$};

\draw[black,thick] (3.5,1.9) -- (4.6,1.9);
\draw[black,thick] (3.5,2.4) -- (4.6,2.4);
\draw (3.5,1.9) node[left] {$3$};
\draw (3.5,2.4) node[left] {$3$};

\draw (5.0,1.375) node {\Large $\mapsto$};

\draw[black] (5.4,0) rectangle (6.4,1);

\foreach \j in {0,...,3}{
  \draw[black] (5.4+0.25*\j,1.65) rectangle (5.4+0.25*\j+0.25,2.65);
}
\draw[black,thick] (5.4+0.125*4,1.65) -- (5.4+0.125*4,2.85);
\draw (5.4+0.125*4,2.85) node[above] {$0$};

\draw[black,thick] (5.4+0.125*2,1.65) -- (5.4+0.125*2,2.8);
\draw[black,thick] (5.4+0.125*6,1.65) -- (5.4+0.125*6,2.8);
\draw (5.4+0.125*2,2.8) node[above] {$2$};
\draw (5.4+0.125*6,2.8) node[above] {$2$};

\draw[black,thick] (5.4+0.125*1,1.65) -- (5.4+0.125*1,2.725);
\draw[black,thick] (5.4+0.125*3,1.65) -- (5.4+0.125*3,2.725);
\draw[black,thick] (5.4+0.125*5,1.65) -- (5.4+0.125*5,2.725);
\draw[black,thick] (5.4+0.125*7,1.65) -- (5.4+0.125*7,2.725);
\draw (5.4+0.125*1,2.725) node[above] {$4$};
\draw (5.4+0.125*3,2.725) node[above] {$4$};
\draw (5.4+0.125*5,2.725) node[above] {$4$};
\draw (5.4+0.125*7,2.725) node[above] {$4$};

\draw[black,thick] (5.2,2.15) -- (6.4,2.15);
\draw (5.2,2.15) node[left] {$1$};

\draw[black,thick] (5.25,1.9) -- (6.4,1.9);
\draw[black,thick] (5.25,2.4) -- (6.4,2.4);
\draw (5.25,1.9) node[left] {$3$};
\draw (5.25,2.4) node[left] {$3$};

\draw[black,thick] (5.325,1.65+1*0.125) -- (6.4,1.65+1*0.125);
\draw[black,thick] (5.325,1.65+3*0.125) -- (6.4,1.65+3*0.125);
\draw[black,thick] (5.325,1.65+5*0.125) -- (6.4,1.65+5*0.125);
\draw[black,thick] (5.325,1.65+7*0.125) -- (6.4,1.65+7*0.125);
\draw[black,thick] (5.325,1.65+1*0.125) node[left] {$5$};
\draw[black,thick] (5.325,1.65+3*0.125) node[left] {$5$};
\draw[black,thick] (5.325,1.65+5*0.125) node[left] {$5$};
\draw[black,thick] (5.325,1.65+7*0.125) node[left] {$5$};

\end{tikzpicture}
\end{center}
\caption{The data distributions of the source domain (bottom) and target domain
(top) throughout the execution of a 2D butterfly algorithm with $p=N^2=64$ 
expressed using bitwise process rank partitions. Notice that the product space 
remains evenly distributed throughout the entire computation.}
\label{fig:2D-par-bin}
\end{figure}
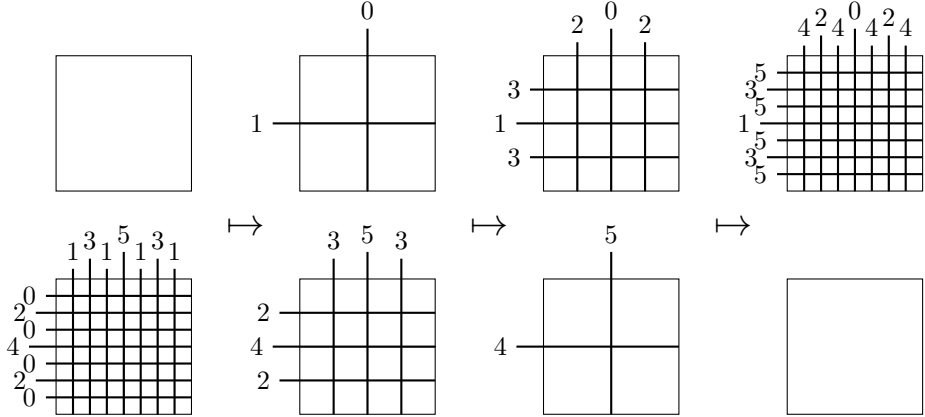

The bit-reversal viewpoint is especially useful when discussing the parallel 
algorithm in higher dimensions, as more bookkeeping is required in order to 
precisely describe distribution of the product space of two multi-dimensional 
domains.
It can be seen from the one-dimensional case shown in 
Fig.~\ref{fig:1D-par-bin} that the bitwise rank-reversal is the result of 
each stage of the parallel algorithm moving the finest-scale bitwise partition 
of the source domain onto the target domain.
In order to best visualize the structure of the equivalent process in higher 
dimensions, it is useful to switch from the matrix-centric viewpoint used in 
the 1D example of Fig.~\ref{fig:1D-par-bin} to the dyadic viewpoint of the 
two-dimensional example of Fig.~\ref{fig:2D-par-bin}.

The main generalization required for the multi-dimensional algorithm is that, 
rather than pairs of processes interacting, $2^d$ processes will need to work 
together in order to map $2^d$ sets of weights from one level to the next.
Just as in the one-dimensional case, each process need only receive one of the 
$2^d$ sets of weights of the result, and the appropriate generalization of each
pairwise exchange is a call to \verb!MPI_Reduce_scatter_block! over a team of 
$2^d$ processes, which only requires each process to send $d$ 
messages and $(2^d-1) r$ entries of data, and to perform $(2^d-1) r$ 
flops~\cite{TRG-coll}.
This communication pattern is precisely the mechanism in which the bitwise 
partitions are moved from the source domain to the target domain: teams of 
processes owning weights which neighbor in the source domain cooperate in order
to produce the weights needed for the next level, which neighbor in the target 
domain.

It is thus easy to see that the per-process cost of the multi-dimensional 
parallel butterfly algorithm is also at most 
$O((r^2 + \beta r + \alpha) \log N)$.
When analytical tensor-product interpolation is used, the local computation 
required for linearly mapping $2^d$ sets of weights from one stage to the next 
can be as low as $O(r^{1+1/d})$, which clearly lowers the parallel complexity 
to $O((r^{1+1/d} + \beta r + \alpha) \log N)$. 
However, there is an additional $O(r^2)$ cost for transitioning from 
column-space to row-space interpolation in the middle of the algorithm when 
the phase function does not have tensor-product structure~\cite{CDY-butterfly}.

\subsection{Algorithm for $p=N^d$ processes}
In order to give precise pseudocode, it is useful to codify the bitwise 
bisections applied to $X$ and $Y$ in terms of two {\em stacks}, $\mathcal{D}_X$
and $\mathcal{D}_Y$.
Each item in the stack is then uniquely specified by the dimension of the domain
it is being applied to and the bit of the process rank used to determine which
half of the domain it will be moved to.
For instance, to initialize the bisection stack $\mathcal{D}_Y$ for the
two-dimensional $p=64$ case shown in Fig.~\ref{fig:2D-par-bin}, we might 
run the following steps:

\begin{algorithm}[H]
\DontPrintSemicolon
$\mathcal{D}_X := \mathcal{D}_Y := \varnothing$\;
\For{$j=0,...,\log_2 p - 1$}{
  $\mathcal{D}_Y.\text{push}((j \bmod d,(\log_2 p-1)-j))$\;
}
\caption{Initialize bisection stacks in $d$-dimensions with $p$ processes}
\label{alg:bisect}
\end{algorithm}
\noindent
This algorithm effectively cycles through the $d$ dimensions bisecting based
upon each of the $\log_2 p$ bits of the process ranks, starting with the 
most-significant bit (with index $\log_2 p - 1$).

Once the bisection stacks have been initialized, the data distribution at level
$\ell$ is defined by sequentially {\em popping} the $d \ell$ bisections off the
top of the $\mathcal{D}_Y$ stack and {\em pushing} them onto the $\mathcal{D}_X$
stack, which we may express as running 
$\mathcal{D}_X.\text{push}(\mathcal{D}_Y.\text{pop}())$ $d \ell$ times.
We also make use of the notation $\mathcal{D}_X(q)$ and $\mathcal{D}_Y(q)$ for 
the portions of $X$ and $Y$ which $\mathcal{D}_X$ and $\mathcal{D}_Y$ 
respectively assign to process $q$.
Lastly, given $p$ processes, we define the {\em bit-masking operator}, 
$\mathcal{M}_a^b(q)$, as
\begin{equation}
  \mathcal{M}_a^b(q) = 
  \{ n \in [0,p) : \text{bit}_j(n)=\text{bit}_j(q),\;\forall j \not\in [a,b) \},
\end{equation}
so that, as long as $0 \le a \le b < \log_2 p$, the cardinality of 
$\mathcal{M}_a^b(q)$, denoted $\text{card}(\mathcal{M}_a^b(q))$, is $2^{b-a}$.
In fact, it can be seen that stage $\ell$ of the $p=N^d$ parallel butterfly 
algorithm, listed as Algorithm~\ref{alg:parallel-equal}, requires process $q$ 
to perform a reduce-scatter summation over the team of $2^d$ processes denoted 
by $\mathcal{M}_{d \ell}^{d(\ell+1)}(q)$.

\begin{algorithm}
\DontPrintSemicolon
\tcp{Initialize bitwise-bisection stacks}
$\mathcal{D}_X := \mathcal{D}_Y := \varnothing$\;
\For{$j=0,...,\log_2 p - 1$}{
  $\mathcal{D}_Y.\text{push}(((j \bmod d,(\log_2 p-1)-j))$\;
}
$\left( \mathcal{D}_X(q) = X,\; \cup_q \mathcal{D}_Y(q) = Y \right)$\;
$A := \mathcal{D}_X(q),\; B := \mathcal{D}_Y(q)$\;
\tcp{Initialize local weights}
$\mathsf{w_{XB} := V_{XB}\, g_B}$\;
\tcp{Translate the weights in parallel}
\For{$\ell=0,...,\log_2(N)-1$}{
  $\{A_c\}_{c=0}^{2^d-1} := \text{children}(A),\;\; B_p := \text{parent}(B)$\;
  \ForEach{$c=0,...,2^d-1$}{
    $\mathsf{w_{A_c B_p} := T_{A_c B_p}^{AB} w_{AB}}$\;
  }
  \For{$j=0,...,d-1$}{
    $\mathcal{D}_X.\text{push}(\mathcal{D}_Y.\text{pop}())$\;
  }
  $A := \mathcal{D}_X(q),\; B := \mathcal{D}_Y(q)$\;
  $\mathsf{w_{AB}} := 
    \text{SumScatter}(\{\mathsf{w_{A_c B_p}}\}_{c=0}^{2^d-1},
                      \mathcal{M}_{d\ell}^{d(\ell+1)}(q))$\;
}
$\left(\cup_q \mathcal{D}_X(q) = X,\; \mathcal{D}_Y(q) = Y\right)$\;
\tcp{Final evaluation}
$\mathsf{f_{AY} := U_{AY}\, w_{AY}}$\;
\caption{Parallel butterfly algorithm with $p=N^d$ process from the point of 
 view of process $q$.}
\label{alg:parallel-equal}
\end{algorithm}

\subsection{Algorithm for $p \le N^d$ processes}
We will now generalize the algorithm of the previous subsection to any case
where a power-of-two number of processes less than or equal to $N^d$ is used
and show that the cost is at most
\[
  O\left(r^2 \frac{N^d}{p} \log N + 
   \left(\beta r\frac{N^d}{p}+\alpha\right)\log p\right),
\]
where $N^d/p$ is the number of interactions assigned to each process during
each stage of the algorithm.
This cost clearly reduces to that of the previous subsection when $p=N^d$ and
to the $O(r^2 N^d \log N)$ cost of the sequential algorithm when $p=1$.

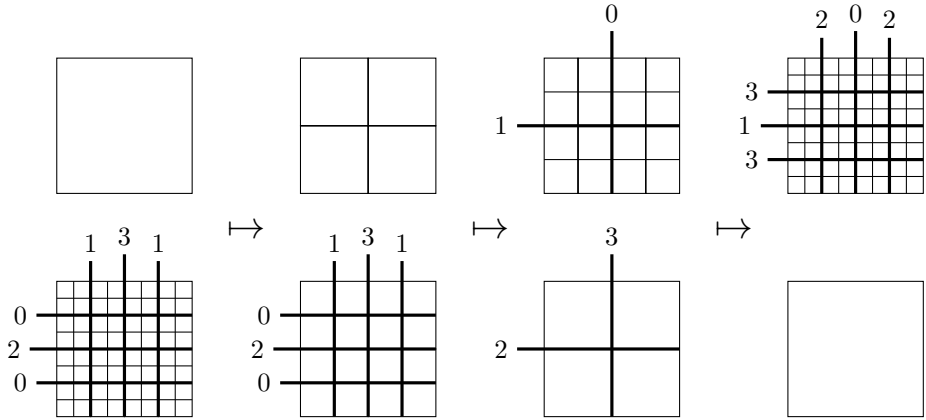
\begin{figure}
\begin{center}
\begin{tikzpicture}[scale=1.8]

\foreach \j in {0,...,3}{
  \draw[black] (0.25*\j,0) rectangle (0.25*\j+0.25,1);
}
\draw[black,very thick] (0.125*4,0) -- (0.125*4,1.2);
\draw (0.125*4,1.2) node[above] {$3$};

\draw[black,very thick] (0.125*2,0) -- (0.125*2,1.15);
\draw[black,very thick] (0.125*6,0) -- (0.125*6,1.15);
\draw (0.125*2,1.15) node[above] {$1$};
\draw (0.125*6,1.15) node[above] {$1$};

\draw (0.125*1,0) -- (0.125*1,1);
\draw (0.125*3,0) -- (0.125*3,1);
\draw (0.125*5,0) -- (0.125*5,1);
\draw (0.125*7,0) -- (0.125*7,1);

\draw[black,very thick] (-0.2,0.5) -- (1,0.5);
\draw (-0.2,0.5) node[left] {$2$};

\draw[black,very thick] (-0.15,0.25) -- (1,0.25);
\draw[black,very thick] (-0.15,0.75) -- (1,0.75);
\draw (-0.15,0.25) node[left] {$0$};
\draw (-0.15,0.75) node[left] {$0$};

\draw (0,1*0.125) -- (1,1*0.125);
\draw (0,3*0.125) -- (1,3*0.125);
\draw (0,5*0.125) -- (1,5*0.125);
\draw (0,7*0.125) -- (1,7*0.125);

\draw[black] (0,1.65) rectangle (1,2.65);

\draw (1.4,1.375) node {\Large $\mapsto$};

\draw[black] (1.8+0.5*0,0.5*0) rectangle (1.8+0.5*1,0.5*0+0.5);
\draw[black] (1.8+0.5*1,0.5*0) rectangle (1.8+0.5*2,0.5*0+0.5);
\draw[black] (1.8+0.5*0,0.5*1) rectangle (1.8+0.5*1,0.5*1+0.5);
\draw[black] (1.8+0.5*1,0.5*1) rectangle (1.8+0.5*2,0.5*1+0.5);

\draw[black,very thick] (1.8+0.125*4,0) -- (1.8+0.125*4,1.2);
\draw (1.8+0.125*4,1.2) node[above] {$3$};

\draw[black,very thick] (1.8+0.125*2,0) -- (1.8+0.125*2,1.15);
\draw[black,very thick] (1.8+0.125*6,0) -- (1.8+0.125*6,1.15);
\draw (1.8+0.125*2,1.15) node[above] {$1$};
\draw (1.8+0.125*6,1.15) node[above] {$1$};

\draw[black,very thick] (1.6,0.5) -- (2.8,0.5);
\draw (1.6,0.5) node[left] {$2$};

\draw[black,very thick] (1.65,0.25) -- (2.8,0.25);
\draw[black,very thick] (1.65,0.75) -- (2.8,0.75);
\draw (1.65,0.25) node[left] {$0$};
\draw (1.65,0.75) node[left] {$0$};

\draw[black] (1.8+0.5*0,1.65+0.5*0) rectangle (1.8+0.5*1,0.5*0+2.15);
\draw[black] (1.8+0.5*1,1.65+0.5*0) rectangle (1.8+0.5*2,0.5*0+2.15);
\draw[black] (1.8+0.5*0,1.65+0.5*1) rectangle (1.8+0.5*1,0.5*1+2.15);
\draw[black] (1.8+0.5*1,1.65+0.5*1) rectangle (1.8+0.5*2,0.5*1+2.15);

\draw (1.8+0.125*4,1.65) -- (1.8+0.125*4,2.65);

\draw (1.8,2.15) -- (2.8,2.15);

\draw (3.2,1.375) node {\Large $\mapsto$};

\draw[black] (3.6,0) rectangle (4.6,1);

\draw[black,very thick] (3.6+0.125*4,0) -- (3.6+0.125*4,1.2);
\draw (3.6+0.125*4,1.2) node[above] {$3$};

\draw[black,very thick] (3.4,0.5) -- (4.6,0.5);
\draw (3.4,0.5) node[left] {$2$};

\foreach \j in {0,...,3}{
  \draw[black] (3.6+0.25*\j,1.65) rectangle (3.6+0.25*\j+0.25,2.65);
}
\draw[black,very thick] (3.6+0.125*4,1.65) -- (3.6+0.125*4,2.85);
\draw (3.6+0.125*4,2.85) node[above] {$0$};

\draw (3.6+0.125*2,1.65) -- (3.6+0.125*2,2.65);
\draw (3.6+0.125*6,1.65) -- (3.6+0.125*6,2.65);

\draw[black,very thick] (3.4,2.15) -- (4.6,2.15);
\draw (3.4,2.15) node[left] {$1$};

\draw (3.6,1.9) -- (4.6,1.9);
\draw (3.6,2.4) -- (4.6,2.4);

\draw (5.0,1.375) node {\Large $\mapsto$};

\draw[black] (5.4,0) rectangle (6.4,1);

\foreach \j in {0,...,3}{
  \draw[black] (5.4+0.25*\j,1.65) rectangle (5.4+0.25*\j+0.25,2.65);
}
\draw[black,very thick] (5.4+0.125*4,1.65) -- (5.4+0.125*4,2.85);
\draw (5.4+0.125*4,2.85) node[above] {$0$};

\draw[black,very thick] (5.4+0.125*2,1.65) -- (5.4+0.125*2,2.8);
\draw[black,very thick] (5.4+0.125*6,1.65) -- (5.4+0.125*6,2.8);
\draw (5.4+0.125*2,2.8) node[above] {$2$};
\draw (5.4+0.125*6,2.8) node[above] {$2$};

\draw (5.4+0.125*1,1.65) -- (5.4+0.125*1,2.65);
\draw (5.4+0.125*3,1.65) -- (5.4+0.125*3,2.65);
\draw (5.4+0.125*5,1.65) -- (5.4+0.125*5,2.65);
\draw (5.4+0.125*7,1.65) -- (5.4+0.125*7,2.65);

\draw[black,very thick] (5.2,2.15) -- (6.4,2.15);
\draw (5.2,2.15) node[left] {$1$};

\draw[black,very thick] (5.25,1.9) -- (6.4,1.9);
\draw[black,very thick] (5.25,2.4) -- (6.4,2.4);
\draw (5.25,1.9) node[left] {$3$};
\draw (5.25,2.4) node[left] {$3$};

\draw (5.4,1.65+1*0.125) -- (6.4,1.65+1*0.125);
\draw (5.4,1.65+3*0.125) -- (6.4,1.65+3*0.125);
\draw (5.4,1.65+5*0.125) -- (6.4,1.65+5*0.125);
\draw (5.4,1.65+7*0.125) -- (6.4,1.65+7*0.125);

\end{tikzpicture}
\end{center}
\caption{A 2D parallel butterfly algorithm with $N^d=64$ and $p=16$. The first 
stage does not require any communication, and the remaining 
$\lceil\log_{2^d} p\rceil=2$ 
stages each require teams of $2^d$ processes to coordinate in order to perform
$N^d/p=4$ simultaneous linear transformations.}
\label{fig:2D-par-bin-less}
\end{figure}

Consider the case where $p < N^d$, such as the 2D example in 
Fig.~\ref{fig:2D-par-bin-less}, where $p=16$ and $N^2=64$.
As shown in the bottom-left portion of the figure, each process is initially 
assigned a contiguous region of the source domain which contains multiple 
source boxes.
If $p$ is significantly less than $N^d$, it can be seen that the first several
stages of the butterfly algorithm will not require any communication, as 
the linear transformations which coarsen the source domain (and refine the 
target domain) will take place over sets of source boxes assigned to a single 
process.
More specifically, when $N$ is a power of two, we may decompose the number of 
stages of the butterfly algorithm, $\log_2 N$, as
\[
  \log_2 N 
  = 
  \lfloor \log_{2^d} \frac{N^d}{p} \rfloor + 
  \lceil \log_{2^d} p \rceil,
\]
where the former term represents the number of stages at the beginning of the 
butterfly algorithm which may be executed entirely locally, and the latter 
represents the number of stages requiring communication.
In the case of Fig.~\ref{fig:2D-par-bin-less}, $\log_{2^d}(N^d/p)=1$, 
and only the first stage does not require communication.

It is important to notice that each of the last $\lceil \log_{2^d} p \rceil$
stages still requires at most $2^d$ processes to interact, and so the only 
difference in the communication cost is that $O(r N^d/p)$ entries of data need
to be exchanged within the team rather than just $O(r)$.
Note that any power of two number of processes can be used with such an 
approach, though, when $p$ is not an integer power of $2^d$, the first 
communication step will involve between $2$ and $2^{d-1}$ processes, while the 
remaining stages will involve teams of $2^d$ processes.
Despite these corner cases, the total communication cost can easily be seen to 
be $O((\beta r (N^d/p) + \alpha)\log p)$, and the total computation cost is
at most $O(r^2 (N^d/p) \log N)$.

Algorithm~\ref{alg:parallel-less} gives a precise prescription of these ideas
using the bisection stacks and bit-masking operators defined in the previous 
subsection.
Because each process can now have multiple source and target boxes 
assigned to it at each stage of the algorithm, we denote these sets of boxes
as
\begin{eqnarray*}
  \mathcal{T}_X(\ell)|_{\mathcal{D}_X(q)} &=& 
  \{ A \in \mathcal{T}_X(\ell) : A \in \mathcal{D}_X(q) \},\,\text{and} \\
  \mathcal{T}_Y(\ell)|_{\mathcal{D}_Y(q)} &=& 
  \{ A \in \mathcal{T}_Y(\ell) : A \in \mathcal{D}_Y(q) \}.
\end{eqnarray*}

\begin{algorithm}
\DontPrintSemicolon
\tcp{Initialize bitwise-bisection stacks}
$\mathcal{D}_X := \mathcal{D}_Y := \varnothing$\;
\For{$j=0,...,\log_2 p - 1$}{
  $\mathcal{D}_Y.\text{push}(((j \bmod d,(\log_2 p-1)-j))$\;
}
$\left( \mathcal{D}_X(q) = X,\;
 \cup_q \mathcal{D}_Y(q) = Y,\; 
 \text{card}\!\left(\mathcal{T}_Y(0)|_{\mathcal{D}_Y(q)}\right)=\frac{N^d}{p} 
 \right)$\;
$\mathcal{A} := \mathcal{T}_X(0)|_{\mathcal{D}_X(q)},\; 
 \mathcal{B} := \mathcal{T}_Y(0)|_{\mathcal{D}_Y(q)}$\;
\tcp{Initialize local weights}
\ForEach{$B \in \mathcal{B}$}{
  $\mathsf{w_{XB} := V_{XB}\, g_B}$\;
}
\tcp{Translate weights in parallel}
$s := \log_2\left( \frac{N^d}{p} \bmod 2^d \right)$\;
\For{$\ell=0,...,\log_2 N - 1$}{
  $\tilde{\mathcal{A}} := \text{children}(\mathcal{A}),\;
   \tilde{\mathcal{B}} := \text{parents}(\mathcal{B})$\;
  \ForEach{$(\tilde A,\tilde B) \in 
            \tilde{\mathcal{A}} \times \tilde{\mathcal{B}}$}{
    $\mathsf{w_{\tilde A \tilde B} := 0}$\;
  }
  \ForEach{$(A,B) \in \mathcal{A} \times \mathcal{B}$}{
    $\{A_c\}_{c=0}^{2^d-1} := \text{children}(A),\; B_p := \text{parent}(B)$\;
    \ForEach{$c=0,...,2^d-1$}{
      $\mathsf{w_{A_c B_p}\, +\!\!= T_{A_c B_p}^{A B} w_{AB}}$\;
    }
  }
  \If{$\ell < \lfloor \log_{2^d} \frac{N^d}{p} \rfloor$}{
    $\mathcal{A} := \tilde{\mathcal{A}},\; \mathcal{B} := \tilde{\mathcal{B}}$\;
  }
  \ElseIf{$\ell = \lfloor \log_{2^d} \frac{N^d}{p} \rfloor\, 
    \text{and}\, s \neq 0$}{
    \For{$j=0,...,s-1$}{
      $\mathcal{D}_X.\text{push}(\mathcal{D}_Y.\text{pop}())$\;
    }
    $\mathcal{A} := \mathcal{T}_X(\ell+1)|_{\mathcal{D}_X(q)},\; 
     \mathcal{B} := \mathcal{T}_Y(\ell+1)|_{\mathcal{D}_Y(q)}$\;
    $\{\mathsf{w_{AB}}\}_{A\in\mathcal{A},B\in\mathcal{B}} := 
     \text{SumScatter}(\{\mathsf{w_{\tilde A\tilde B}}\}_{
       \tilde A \in \tilde{\mathcal{A}},\tilde B \in \tilde{\mathcal{B}}},
       \mathcal{M}_{d\ell}^{d\ell+s}(q))$\;
  }
  \Else{
    \For{$j=0,...,d-1$}{
      $\mathcal{D}_X.\text{push}(\mathcal{D}_Y.\text{pop}())$\;
    }
    $\mathcal{A} := \mathcal{T}_X(\ell+1)|_{\mathcal{D}_X(q)},\; 
     \mathcal{B} := \mathcal{T}_Y(\ell+1)|_{\mathcal{D}_Y(q)}$\;
    $\{\mathsf{w_{AB}}\}_{A\in\mathcal{A},B\in\mathcal{B}} := 
     \text{SumScatter}(\{\mathsf{w_{\tilde A\tilde B}}\}_{
       \tilde A \in \tilde{\mathcal{A}},\tilde B \in \tilde{\mathcal{B}}},
       \mathcal{M}_{d(\ell-1)+s}^{d\ell+s}(q))$\;
  }
}
$\left(\text{card}\!\left(\mathcal{T}_X(\log_2 N)|_{\mathcal{D}_X(q)}\right)=
       \frac{N^d}{p},\; 
  \cup_q \mathcal{D}_X(q) = X,\; \mathcal{D}_Y(q) = Y\right)$\;
\tcp{Final evaluations}
\ForEach{$A \in \mathcal{A}$}{
  $\mathsf{f_{AY} := U_{AY}\, w_{AY}}$\;
}
\caption{Parallel butterfly algorithm for $p \le N^d$ processes from the view of
 process $q$. $p$ is assumed to be a power of two.}
\label{alg:parallel-less}
\end{algorithm}


\section{Experimental results}
As previously mentioned, our performance experiments focus on the class of 
integral operators whose kernels are of the form of Eq.~\eqref{phase}.
While this may seem overly restrictive, a large class of important transforms 
falls into this category, most notably:
the Fourier transform, where $\Phi(x,y) = 2 \pi x \cdot y$, 
backprojection~\cite{Demanet-SAR}, 
hyperbolic Radon transforms~\cite{Hu-generalized}, 
and Egorov operators, which then provide a means of efficiently applying 
Fourier Integral Operators~\cite{CDY-butterfly}.
Due to the extremely special (and equally delicate) structure of Fourier 
transforms, a number of highly-efficient parallel algorithms already 
exist for both uniform~\cite{FosterWorley-STM,Pippig-PFFT,Czechowki-3DFFT} and 
non-uniform~\cite{PippigPotts-PNFFT} Fourier transforms, and so we will instead 
concentrate on more sophisticated kernels.
We note that the high-level communication pattern and costs of the parallel 1D FFT 
mentioned in \cite{FosterWorley-STM} are closely related to those of our parallel 1D 
butterfly algorithm.

Algorithm~\ref{alg:parallel-less} was instantiated in the new 
\verb!DistButterfly! library using black-box, user-defined phase functions, 
and the low-rank approximations and translation operators introduced in 
\cite{CDY-butterfly}.
The library was written using \verb!C++11! in order to template the 
implementation over the dimension of the problem, and all inter-process 
communication was expressed via the Message Passing Interface (MPI).
All tests were performed on the Argonne Leadership Computing Facility's 
Blue Gene/Q installation using a port of Clang and IBM's ESSL 5.1, and kernel 
evaluations were accelerated by batching them together and calling 
ESSL's MASS routines, \verb!vsin!, \verb!vcos!, and \verb!vsincos!.
All calculations were performed with (64-bit) double-precision arithmetic.

The current implementation is written purely with MPI and does not employ 
multi-threading, and because Blue Gene/Q's cores require multiple threads in 
order to saturate the machine's floating-point units~\cite{BLIS-TR}, 
it was empirically found that launching four MPI processes for each of the 16 
cores on each node resulted in the best performance.
All strong-scaling tests were therefore conducted in the range of 
1-node/16-cores/64-processes and 1024-nodes/16,384-cores/65,536-processes.

\subsection{Hyperbolic Radon transforms}
Our first set of experiments used a phased function of the form
\[
  \Phi((x_0,x_1),(h,p))=
  2\pi p \sqrt{x_0^2 + x_1^2 h^2},
\]
which corresponds to a 2D (or 3D)\footnote{If the degrees of freedom in the 
second and third dimensions are combined~\cite{Hu-generalized}.} 
{\em hyperbolic Radon transform}~\cite{Hu-generalized}
and has many applications in seismic imaging.
For the strong-scaling tests, the numerical rank was fixed at $r=4^2$, which
corresponds to a tensor product of two four-point Chebyshev grids for each 
low-rank interaction.

As can be seen from Fig.~\ref{fig:perf-hypradon}, the performance of both the 
$N=128$ and $N=256$ problems continued to improve all the way to the $p=N^2$
limits (respectively, 16,384 and 65,536 processes).
As expected, the larger problems display the best strong scaling, and the 
largest problem which would fit on one node, $N=1024$, scaled from 64 to 
65,536 processes with roughly 90.5\% efficiency.

\begin{figure}
\centering
\begin{tikzpicture}
\begin{loglogaxis}[width=10cm,
  xlabel={Number of cores of Blue Gene/Q},
  ylabel={Walltime [seconds]},
  xmin=16,xmax=16384,ymin=0.002,ymax=300,clip=false,
  cycle list name=journal]

%
%

\addplot coordinates {
  (16,0.6745)
  (32,0.3544)
  (64,0.1798)
  (128,0.0925)
  (256,0.0474)
  (512,0.0245)
  (1024,0.0129)
  (2048,0.0073)
  (4096,0.0046)
};
\addplot coordinates {
  (16,0.6745)
  (32,0.3372)
  (64,0.1686)
  (128,0.0843)
  (256,0.0422)
  (512,0.0211)
  (1024,0.0105)
  (2048,0.0053)
  (4096,0.0026)
};

%
%

\addplot coordinates {
  (16,3.0036)
  (32,1.5115)
  (64,0.7744)
  (128,0.4054)
  (256,0.2031)
  (512,0.1063)
  (1024,0.0540)
  (2048,0.0284)
  (4096,0.0151)
  (8192,0.0085)
  (16384,0.0049)
};
\addplot coordinates {
  (16,3.0036)
  (32,1.5018)
  (64,0.7509)
  (128,0.3755)
  (256,0.1877)
  (512,0.0939)
  (1024,0.0469)
  (2048,0.0235)
  (4096,0.0117)
  (8192,0.0059)
  (16384,0.0029)
};

%
%

\addplot coordinates {
  (16,13.8999)
  (32,6.9219)
  (64,3.4674)
  (128,1.7500)
  (256,0.8937)
  (512,0.4667)
  (1024,0.2364)
  (2048,0.1211)
  (4096,0.0618)
  (8192,0.0322)
  (16384,0.0169)
};
\addplot coordinates {
  (16,13.8999)
  (32,6.9500)
  (64,3.4750)
  (128,1.7375)
  (256,0.8687)
  (512,0.4344)
  (1024,0.2172)
  (2048,0.1086)
  (4096,0.0543)
  (8192,0.0271)
  (16384,0.0136)
};

%
%

\addplot coordinates {
  (16,63.7311)
  (32,30.6760)
  (64,15.7723)
  (128,7.7600)
  (256,3.9071)
  (512,1.9542)
  (1024,1.0028)
  (2048,0.5245)
  (4096,0.2618)
  (8192,0.1364)
  (16384,0.0688)
};
\addplot coordinates {
  (16,63.7311)
  (32,31.8655)
  (64,15.9328)
  (128,7.9664)
  (256,3.9832)
  (512,1.9916)
  (1024,0.9958)
  (2048,0.4979)
  (4096,0.2489)
  (8192,0.1245)
  (16384,0.0622)
};

%
%

\addplot coordinates {
  (64,72.0656)
  (128,34.2627)
  (256,17.8815)
  (512,8.5020)
  (1024,4.3126)
  (2048,2.1954)
  (4096,1.1002)
  (8192,0.5710)
  (16384,0.3046)
};
\addplot coordinates {
  (64,72.0656)
  (128,36.0328)
  (256,18.0164)
  (512,9.0082)
  (1024,4.5041)
  (2048,2.2521)
  (4096,1.1260)
  (8192,0.5630)
  (16384,0.2815)
};

%
%

\addplot coordinates {
  (256,79.1480)
  (512,36.9299)
  (1024,18.6219)
  (2048,9.6551)
  (4096,4.7415)
  (8192,2.3768)
  (16384,1.2201)
};
\addplot coordinates {
  (256,79.1480)
  (512,39.5740)
  (1024,19.7870)
  (2048,9.8935)
  (4096,4.9467)
  (8192,2.4734)
  (16384,1.2367)
};

%
%

\addplot coordinates {
  (1024,80.1919)
  (2048,41.0508)
  (4096,21.0158)
  (8192,10.2855)
  (16384,5.1880)
};
\addplot coordinates {
  (1024,80.1919)
  (2048,40.0960)
  (4096,20.0480)
  (8192,10.0240)
  (16384,5.0120)
};

%
%

\addplot coordinates {
  (4096,93.4758)
  (8192,46.8329)
  (16384,22.1774)
};
\addplot coordinates {
  (4096,93.4758)
  (8192,46.7379)
  (16384,23.3690)
};

\addplot coordinates {
  (16384,101.4000)
};
\end{loglogaxis}
\end{tikzpicture}
\caption{A combined strong-scaling and asymptotic complexity test for 
Algorithm~\ref{alg:parallel-less} using analytical interpolation for a 
2D hyperbolic Radon transform with numerical ranks of $r=4^2$.
From bottom-left to top-right, the tests involved $N^2$ source and 
target boxes with $N$ equal to 128, 256, ..., 32768. Note that four MPI 
processes were launched per core in order to maximize performance and that each dashed line 
corresponds to linear scaling relative to the test using the smallest possible portion of 
the machine.}

\label{fig:perf-hypradon}
\end{figure}

%

\subsection{3D Generalized Radon analogue}
Our second test involves an analogue of a 3D generalized Radon 
transform~\cite{Beylkin-genRadon} and again has applications in seismic imaging.
As in \cite{CDY-butterfly}, we use a phase function of the form
\[
  \Phi(x,p) = \pi \left(x \cdot p + 
    \sqrt{\gamma(x,p)^2 + \kappa(x,p)^2}\right),
\]
where
\begin{eqnarray*}
 \gamma(x,p) &=& p_0 \left(2+\sin(2\pi x_0) \sin(2\pi x_1)\right)/3,
 \,\text{and}\\
 \kappa(x,p) &=& p_1 \left(2+\cos(2\pi x_0) \cos(2\pi x_1)\right)/3.
\end{eqnarray*}
Notice that, if not for the non-linear square-root term, 
the phase function would be equivalent to that of a Fourier transform.
The numerical rank was chosen to be $r=5^3$ in order to provide roughly one 
percent relative error in the supremum norm, and due to both the increased rank,
increased dimension, and more expensive phase function, the runtimes are 
significantly higher than those of the previous example for cases with 
equivalent numbers of source and target boxes.

Just as in the previous example, the smallest two problem sizes, $N=16$ and 
$N=32$, were observed to strong scale up to the limit of $p=N^3$ (respectively,
4096 and 32,768) processes. 
The most interesting problem size is $N=64$, which corresponds to the largest
problem which was able to fit in the memory of a single node.
The strong scaling efficiency from one to 1024 nodes was observed to be 82.3\%
in this case.

\begin{figure}
\centering
\begin{tikzpicture}
\begin{loglogaxis}[width=10cm,
  xlabel={Number of cores of Blue Gene/Q},
  ylabel={Walltime [seconds]},
  xmin=16,xmax=16384,ymin=0.02,ymax=300,clip=false,
  cycle list name=journal]

%
%

\addplot coordinates {
  (16,1.0924)
  (32,0.5601)
  (64,0.2824)
  (128,0.1487)
  (256,0.0791)
  (512,0.0428)
  (1024,0.0249)
};
\addplot coordinates {
  (16,1.0924)
  (32,0.5462)
  (64,0.2731)
  (128,0.1366)
  (256,0.0683)
  (512,0.0341)
  (1024,0.0171)
};

%
%

\addplot coordinates {
  (16,10.4525)
  (32,5.2882)
  (64,2.7351)
  (128,1.3864)
  (256,0.7100)
  (512,0.3550)
  (1024,0.1841)
  (2048,0.0965)
  (4096,0.0512)
  (8192,0.0296)
};
\addplot coordinates {
  (16,10.4525)
  (32,5.2263)
  (64,2.6131)
  (128,1.3066)
  (256,0.6533)
  (512,0.3266)
  (1024,0.1633)
  (2048,0.0817)
  (4096,0.0408)
  (8192,0.0204)
};

%
%

\addplot coordinates {
  (16,94.4076)
  (32,47.0929)
  (64,23.7254)
  (128,12.1361)
  (256,6.0350)
  (512,3.1708)
  (1024,1.6171)
  (2048,0.8216)
  (4096,0.4037)
  (8192,0.2115)
  (16384,0.1120)
};
\addplot coordinates {
  (16,94.4076)
  (32,47.2038)
  (64,23.6019)
  (128,11.8010)
  (256,5.9005)
  (512,2.9502)
  (1024,1.4751)
  (2048,0.7376)
  (4096,0.3688)
  (8192,0.1844)
  (16384,0.0922)
};

%
%

\addplot coordinates {
  (64,224.0390)
  (128,113.1440)
  (256,56.0694)
  (512,28.2360)
  (1024,14.3938)
  (2048,7.2030)
  (4096,3.6692)
  (8192,1.8925)
  (16384,1.0045)
};
\addplot coordinates {
  (64,224.0390)
  (128,112.0195)
  (256,56.0097)
  (512,28.0049)
  (1024,14.0024)
  (2048,7.0012)
  (4096,3.5006)
  (8192,1.7503)
  (16384,0.8752)
};

%
%

\addplot coordinates {
  (512,251.9040)
  (1024,128.1830)
  (2048,63.1565)
  (4096,31.4749)
  (8192,16.1943)
  (16384,8.4330)
};
\addplot coordinates {
  (512,251.9040)
  (1024,125.9520)
  (2048,62.9760)
  (4096,31.4880)
  (8192,15.7440)
  (16384,7.8720)
};

%
%

\addplot coordinates {
  (4096,283.5680)
  (8192,143.4520)
  (16384,72.6868)
};
\addplot coordinates {
  (4096,283.5680)
  (8192,141.7840)
  (16384,70.8920)
};

\end{loglogaxis}
\end{tikzpicture}
\caption{A combined strong-scaling and asymptotic complexity test for 
Algorithm~\ref{alg:parallel-less} using analytical interpolation for a 3D 
generalized Radon analogue with numerical rank $r=5^3$.
From bottom-left to top-right, the tests involved $N^3$ source and target boxes
with $N$ equal to 16, 32, ..., 512. Note that four MPI processes were launched
per core in order to maximize performance and that each dashed line corresponds to 
linear scaling relative to the test using the smallest possible portion of the machine.}
\label{fig:perf-genradon}
\end{figure}

\section{Conclusions and future work}
A high-level parallelization of both general-purpose and 
analytical-interpolation based butterfly algorithms was presented, with the 
former resulting in a modeled runtime of 
\[
  O\left(r^2 \frac{N^d}{p} \log N + 
  \left(\beta r\frac{N^d}{p}+\alpha\right)\log p\right)
\]
for $d$-dimensional problems with $N^d$ source and target boxes, rank-$r$ 
approximate interactions, and $p \le N^d$ processes 
(with message latency, $\alpha$, and inverse bandwidth, $\beta$).
The key insight is that a careful manipulation of bitwise-partitions of the 
product space of the source and target domains can both keep the data 
(weight vectors) and the computation (weight translations) evenly 
distributed, and teams of at most $2^d$ need interact via a reduce-scatter 
communication pattern at each of the $\log_2(N)$ stages of the algorithm.
Algorithm~\ref{alg:parallel-less} was then implemented in a black-box manner 
for kernels of the form $K(x,y)=\exp(i \Phi(x,y))$, and strong-scaling of 
90.5\% and 82.3\% efficiency from one to 1024 nodes of Blue Gene/Q was 
observed for a 2D hyperbolic Radon transform and 3D generalized Radon analogue,
respectively.

While, to the best of our knowledge, this is the first parallel implementation 
of the butterfly algorithm, there is still a significant amount of future work.
The most straight-forward of which is to extend the current MPI implementation
to exploit the $N^d/p$-fold trivial parallelism available for the local weight
translations.
This could potentially provide a significant performance improvement on 
``wide'' architectures such as Blue Gene/Q.
Though a significant amount of effort has already been devoted to improving the 
architectural efficiency, e.g., via batch evaluations of phase functions, 
further improvements are almost certainly waiting.

A more challenging direction involves the well-known fact that the butterfly 
algorithm can exploit sparse source and target 
domains~\cite{ChewSong-SFT,Ying-SFT,KunisMelzer-stableButterfly}, and so it 
would be worthwhile to extend our parallel algorithm into this regime.
Lastly, the butterfly algorithm is closely related to the 
{\em directional Fast Multipole Method}~\cite{EngquistYing-oscillatory}, 
which makes use of low-rank interactions over spatial 
cones~\cite{Martinsson-elongated} which can be interpreted as satisfying the 
butterfly condition in angular coordinates. 
It would be interesting to investigate the degree to which our parallelization
of the butterfly algorithm carries over to the directional FMM.

\section*{Availability}
The distributed-memory implementation of the butterfly algorithm for kernels 
of the form $\exp\left(i \Phi(x,y)\right)$, 
\texttt{DistButterfly}, is available at 
\url{github.com/poulson/dist-butterfly} under the GPLv3.
All experiments in this manuscript made use of revision 
\verb!f850f1691c!.
Additionally, parallel implementations of interpolative decompositions are now 
available as part of Elemental~\cite{Poulson-elemental}, which is hosted at 
\url{libelemental.org} under the New BSD License.

\section*{Acknowledgments}
Jack Poulson would like to thank both Jeff Hammond and Hal Finkel for their 
generous help with both \verb!C++11! support and performance issues on 
Blue Gene/Q.
He would also like to thank Mark Tygert for detailed discussions on 
both the history and practical implementation of IDs and related factorizations,
and Gregorio Quintana Ort\'i for discussions on high-performance RRQRs.

This work was partially supported by NSF CAREER Grant No.~0846501 (L.Y.),
DOE Grant No.~DE-SC0009409 (L.Y.), and KAUST. Furthermore, this research used
resources of the Argonne Leadership Computing Facility at Argonne National 
Laboratory, which is supported by the Office of Science of the U.S.~Department
of Energy under contract DE-AC02-06CH11357.

\end{document}